\documentclass[12pt]{article}
\usepackage{amsfonts}

\newtheorem{thm}{Theorem}[section]
\newtheorem{lemma}[thm]{Lemma}
\newtheorem{cor}[thm]{Corollary}

\newtheorem{prop}[thm]{Proposition}
\newtheorem{conjecture}{Conjecture}

\newenvironment{remark}{\par\medskip\noindent{\bf Remark.\ }}{\par\smallskip}
\newcommand{\proof
}{\par\medskip\noindent {\bf Proof.\ \ }}
\newcommand{\ve}[1]{\mathbf{#1}}
\newcommand{\be}{\begin{equation}}
\newcommand{\ee}{\end{equation}}
\newcommand{\openbox}{\leavevmode
  \hbox to8pt{\hfil\vrule\vbox to6pt{\hrule width6pt\vfil\hrule}\vrule}}

\newcommand{\qed}{\hbox to5pt{ } \hfill \openbox\bigskip\medskip}

\newcommand{\cS}{\mbox{$\cal S$}}

\newcommand{\cF}{\mbox{$\cal F$}}
\newcommand{\cG}{\mbox{$\cal G$}}

\newcommand{\cM}{\mbox{$\cal M$}}
\newcommand{\cN}{\mbox{$\cal N$}}

\newcommand{\cE}{\mbox{$\cal E$}}
\newcommand{\Sf}{\mathbb S}

\newcommand{\N}{\mathbb N}

\newcommand{\R}{\mathbb R}

\newcommand{\F}{\mathbb F}

\title{A new proof of a generalization of Gerzon's bound}
\author{G\'abor Heged\"{u}s
\\{\normalsize  \'Obuda University}
\\{\normalsize B\'ecsi \'ut 96, Budapest, Hungary, H-1037}
\\{\normalsize hegedus.gabor@nik.uni-obuda.hu}
}

\begin{document}

\maketitle
\begin{abstract}
In this paper we give a short, new proof of  a natural generalization of Gerzon's bound. This bound improves the Delsarte, Goethals and  Seidel's upper bound in a special case. 

Our proof is a simple application of  the linear algebra bound method.
\end{abstract}
\medskip

\footnotetext{
%{\bf 2010 Mathematics Subject Classification: 52C45, 12Y99, 05D99}
 }

\section{Introduction}

In this paper we give a linear algebraic proof of  the known  upper bound for the size of some special spherical $s$-distance sets. This result generalizes Gerzon's general upper bound for the size of  equiangular spherical set.

%Let $\F$ be a field.
In the following $\R[x_1, \ldots, x_n]=\R[\ve x]$ denotes  the
ring of polynomials in commuting variables $x_1, \ldots, x_n$ over $\R$.

Let  $\cG\subseteq {\R}^n$ be an arbitrary set. 
Denote by $d(\cG)$ the set of (non-zero) distances among points of $\cG$:
$$
d(\cG):=\{d(\ve p_1,\ve p_2):~ \ve p_1,\ve p_2\in \cG,\ve p_1\ne \ve p_2\}.
$$
An $s$-distance set is any subset $\mbox{$\cal H$}\subseteq {\mathbb R}^n$ such that $|d(\mbox{$\cal H$})|\leq s$.

Let $(\ve x,\ve y)$ stand for  the standard scalar product.  Let $s(\cG)$ denote the set of scalar products between the distinct points of $\cG$:
$$
s(\cG):=\{(\ve p_1,\ve p_2):~ \ve p_1,\ve p_2\in \cG,\ve p_1\ne \ve p_2\}.
$$

A {\em spherical} $s$-distance set means any subset $\mbox{$\cal G$}\subseteq {\Sf}^{n-1}$ such that $|s(\mbox{$\cal G$})|\leq s$.

%Bannai, Bannai and  Stanton proved the following result in \cite{BBS} Theorem 1.

%\begin{thm} \label{BBSupper} (Bannai, Bannai and  Stanton \cite{BBS} Theorem 1)
%Suppose that $\cF\subseteq {\R}^n$ is a set satisfying $|d(\cF)|\leq s$. Then
%$$
%|\cF|\leq {n+s\choose s}.
%$$
%\end{thm}
Let $n,s\geq 1$ be integers. Define
$$
M(n,s):={n+s-1\choose s}+{n+s-2\choose s-1}.
$$

Delsarte, Goethals and  Seidel investigated the spherical $s$-distance sets. They  proved a general upper bound for the size of a spherical $s$-distance set in \cite{DGS}.
\begin{thm} (Delsarte, Goethals and  Seidel \cite{DGS})\label{DGGupper}
Suppose that $\cF\subseteq {\Sf}^{n-1}$ is a set satisfying $|s(\cF)|\leq s$.
Then
$$
|\cF|\leq M(n,s).
$$
\end{thm}

Barg and Musin gave an improved upper bpund for the size of a spherical $s$-distance set in a special case in \cite{BM}. Their proof builds upon  Delsarte's ideas (see \cite{DGS}) and they used Gegenbauer polynomials in their argument.
\begin{thm} (Barg, Musin \cite{BM})
Let $n\geq 1$ be a positive integer and  let $s>0$ be an even integer. 
Let $\cS\subseteq {\Sf}^{n-1}$ denote a spherical $s$-distance set  with inner products $t_1,\ldots ,t_s$ such that 
$$
t_1+\ldots +t_s\geq 0.
$$
Then 
$$
|\cS|\leq M(n,s-2)+\frac{n+2s-2}{s}{n+s-1\choose s-1}.
$$
\end{thm}

We point out here the following special case of Theorem \ref{DGGupper}.
\begin{cor} \label{2dist}
Suppose that $\cF\subseteq {\Sf}^{n-1}$ is a set satisfying $|s(\cF)|\leq 2$. Then
$$
|\cS|\leq \frac{n(n+3)}{2}
$$
\end{cor}

An {\em equiangular spherical set} means a two-distance spherical set with scalar products $\alpha$ and $-\alpha$. Let $M(n)$ denote the maximum cardinality of an equiangular spherical set. There is  a very extensive literature devoted to  the determination of  the precise value of $M(n)$ (see \cite{LS}, \cite{GY}). Gerzon gave the first general upper bound for $M(n)$ in \cite{LS} Theorem 8. 

\begin{thm} (Gerzon \cite{LS} Theorem 8)
Let $n\geq 1$ be a positive integer. Then
$$
M(n)\leq \frac{n(n+1)}{2}. 
$$
\end{thm}

Musin proved a stronger version of Gerzon's Theorem in \cite{M} Theorem 1. He used the linear algebra bound method in his proof.
\begin{thm} (Musin, \cite{M} Theorem 1) 
Let $\cS$ be a spherical two-distance set  with inner products $a$ and $b$. Suppose that $a + b\geq  0$. Then
$$
|\cS|\leq \frac{n(n+1)}{2}. 
$$
\end{thm}

De Caen gave a lower bound for the size of an equiangular spherical set in \cite{D}.
\begin{thm} (De Caen \cite{D}) 
Let  $t>0$ be a  positive integer. For each $n = 3 \cdot 2^{2t-1} - 1$ there exists an
equiangular spherical set of $\frac 29 (n + 1)^2$ vectors.
\end{thm}

Our main result is an alternative proof of a natural generalization of Gerzon's bound, which improves the Delsarte, Goethals and  Seidel's upper bound in a special case. Our proof uses the linear algebra bound method. The following statement was proved in \cite{DGS2} Theorem 6.1.  The original proof builds upon   matrix techniques and the addition
formula for Jacobi polynomials.
\begin{thm} (Delsarte, Goethals and  Seidel  \cite{DGS2} Theorem 6.1) \label{main}
Let $s=2t>0$ be an even number and $n>0$ be a positive integer. 
Let $\cS$ denote a spherical $s$-distance set with inner products $a_1,\ldots ,a_t,-a_1,\ldots ,-a_t$. Suppose that $0<a_i< 1$ for each $i$. Then
$$
|\cS|\leq {n+s-1\choose s}.
$$
\end{thm}

%\begin{remark}
%Let $s>0$ be an even integer. Let $n=2$. Then it is easy to verify that the vertices of a regular $s+1$-gon on the unit circle give us an extremal configuration that shows that Theorem \ref{main} is sharp  in this special case. 
%\end{remark}

\section{Preliminaries}

%\subsection{Reduction }

%\subsection{Determinant Criterion }
We prove our main result using the linear algebra bound method and the Determinant Criterion (see \cite{BF} Proposition 2.7). We recall here for the reader's convenience this principle.
\begin{prop} \label{det} (Determinant Criterion)
Let $\F$ denote an arbitrary field. Let $f_i:\Omega \to \F$ be functions and $\ve v_j\in \Omega$ elements for each $1\leq i,j\leq m$  such that the $m \times m$ matrix $B=(f_i(\ve v_j))_{i,j=1}^m$
is non-singular. Then $f_1,\ldots ,f_m$ are linearly independent functions of the space $\F^{\Omega}$. 
\end{prop}

%\subsection{Enumerative combinatorics}

Consider the set of vectors
$$
\cM(n,s):=\{\alpha=(\alpha_1, \ldots , \alpha_n)\in {\N}^n:~ \alpha_1\leq 1, \sum_{i=1}^n \alpha_i \mbox{ is even },\ \sum_{i=1}^n\alpha_i\leq s\}.
$$
Define the set
$$
\cN(n,s)=\{\alpha=(\alpha_1, \ldots , \alpha_n)\in {\N}^n:~ \sum_{i=1}^n \alpha_i\leq s\}.
$$
\begin{lemma} \label{monom}
Let $n,s\geq 1$ be integers.
Then
$$
|\cN(n,s)|={n+s \choose s}.
$$
\end{lemma}
\proof 
For a simple proof of this fact see \cite{CLO} Section 9.2 Lemma 4.\qed

We use the following combinatorial statement in the proofs of our main results.
\begin{cor} \label{count}
Let $n>0$ be a positive integer and $s>0$ be an even integer. 
Then
$$
|\cM(n,s)|={n+s-1\choose s}.
$$
\end{cor}
\proof
It is easy to check there exists a bijection $f:\cM(n,s)\to \cN(n-1,s)$, since $s$ is even.  Namely let $\alpha=(\alpha_1, \ldots , \alpha_n)\in \cM(n,s)$ be an arbitrary element. Then define $f(\alpha):=(\alpha_2, \ldots , \alpha_n)$. It is easy to verify that $f(\cM(n,s))\subseteq \cN(n-1,s)$ and $f$ is a bijection.

Hence $|\cM(n,s)|=|\cN(n-1,s)|$ and Lemma \ref{monom} gives the result.
\qed

\section{The proof}

{\bf Proof of Theorem \ref{main}:}\\   

Consider the real polynomial
$$
g(x_1,\ldots ,x_n)=(\sum_{m=1}^n x_m^2) -1\in \R[x_1, \ldots ,x_n].
$$

Let $\cS=\{\ve v_1, \ldots , \ve v_r\}$ denote a spherical $s$-distance set with inner products $a_1,\ldots ,a_t,-a_1,\ldots ,-a_t$. Here $r=|\cS|$.
Define the polynomials
$$
P_i(x_1, \ldots ,x_n):= \prod_{m=1}^t \Big((\langle \ve x, \ve v_i\rangle)^2-(a_m)^2 \Big)\in \R[\ve x]
$$
for each $1\leq i\leq r$. Clearly $\mbox{deg}(P_i)\leq s=2t$ for each $1\leq i\leq r$. 
%Here $(\ve v_i)_j$ denotes the $j$th coordinate of the vector $\ve v_i$. 

Consider the set of vectors
$$            
\cE(n,s):=\{\alpha=(\alpha_1, \ldots , \alpha_n)\in {\N}^n:~ \sum_i \alpha_i \mbox{ is even },\ \sum \alpha_i\leq s\}
$$
It is easy to verify that if we can expand $P_i$  as a linear combination of monomials, then we get 
%a linear combination of monomials in the form
\begin{equation}  \label{exp} 
P_i(x_1, \ldots ,x_n)=\sum_{\alpha\in \cE(n,s)} c_{\alpha}x^{\alpha},
\end{equation}
where $c_{\alpha}\in \R$ are real coefficients  for each ${\alpha}\in \cE(n,s)$ and $x^{\alpha}$ denotes the monomial  $x_{1}^{\alpha_{1}}\cdot \ldots \cdot x_{n}^{\alpha_{n}}$.  

Since $\ve v_i\in {\Sf}^{n-1}$, this means that the equation
\begin{equation}  \label{sphere} 
x_1^2=1-\sum_{j=2}^n x_j^2
\end{equation}
is true for each $\ve v_i$, where $1\leq i\leq r$. 
%Let $\prec$ an arbitrary term order such that $x_n\prec \ldots \prec x_1$. 
Let $Q_i$ denote the polynomial obtained by writing  $P_i$ as a linear combination of monomials and replacing,
repeatedly, each occurrence of $x_1^t$, where $t\geq 2$,  by a linear combination of other monomials, using the relations (\ref{sphere}).
 
Since $g(\ve v_i)=0$ for each $i$, 
hence $Q_i(\ve v_j)=P_i(\ve v_j)$ for each $1\leq i\ne j\leq r$. 

We prove that the set of polynomials $\{Q_i:~ 1\leq i\leq r\}$ is linearly independent. This fact follows from the Determinant Criterion, when we define $\F:=\R$, $\Omega={\Sf}^{n-1}$ and $f_i:=Q_i$ for each $i$.  It is enough to prove that  $Q_i(\ve v_i)=P_i(\ve v_i)\ne 0$ for each $1\leq i\leq r$ and $Q_i(\ve v_j)=P_i(\ve v_j)=0$ for each $1\leq i\ne j\leq r$, since then we can apply the Determinant Criterion.

But $P_i(\ve v_i)=\prod_{i=1}^m (1-a_m^2)$ and $P_i(\ve v_j)=0$, because $\cS=\{\ve v_1, \ldots , \ve v_r\}$ is a spherical $s$-distance set with inner products $a_1,\ldots ,a_t,-a_1,\ldots ,-a_t$.

%$$
%\cup\{\alpha=(\alpha_1, \ldots , \alpha_n)\in {\N}^n:~ \alpha_1=1, \sum_i \alpha_i \mbox{ is even },\ \sum \alpha_i\leq s\}.
%$$
Then it is easy to check that  we can write $Q_i$ as a linear combination of monomials in the form
$$
Q_i=\sum _{\alpha\in \cM(n,s)} d_{\alpha}x^{\alpha},
$$ 
where $d_{\alpha}\in \R$ are the real coefficients for each ${\alpha}\in \cM(n,s)$. This follows immediately from the expansion (\ref{exp}) and from the relation \ref{sphere}. 
%and ${\rm lm}(g)=x_1^2$.  

Since the polynomials $\{Q_i:~ 1\leq i\leq r\}$ are linearly independent and 
if we expand $Q_i$ as a linear combination of monomials, then  all monomials appearing in this linear combination  contained in the set of monomials 
$$
\{x^{\alpha}:~ \alpha\in \cM(n,s)\}
$$
for each $i$, hence 
$$
r=|\cS|\leq  |\cM(n,s)|={n+s-1\choose s},
$$
by Corollary \ref{count}.
\qed
%\proof
%his is an easy consequence of the following well-known statement  (

\section{Concluding remarks}

The following Conjecture is a natural strengthening of Theorem \ref{main}.
\begin{conjecture} \label{Hconj}
Let $n\geq 1$ be a positive integer and  let $s>0$ be an even integer. 
Let $\cS\subseteq {\Sf}^{n-1}$ denote a spherical $s$-distance set  with inner products $t_1,\ldots ,t_s$ such that 
$$
t_1+\ldots +t_s\geq 0.
$$
Then
$$
|\cS|\leq {n+s-1\choose s}.
$$
\end{conjecture}
%{\bf Acknowledgements.} 

\end{document}